\documentclass[twocolumn,aps,floatfix,superscriptaddress,longbibliography]{revtex4-1}
\pdfoutput=1 
\usepackage{amsmath,amssymb,eucal,graphicx,float,epstopdf,xparse}
\usepackage{epsfig,subfigure}
\usepackage[utf8]{inputenc}

\NewDocumentCommand{\eulerian}{omm}
 {%
  \genfrac<>{0pt}{}{#2}{#3}%
  \IfValueT{#1}{_{\!#1}}%
 }

\usepackage[colorlinks=true, urlcolor=blue, anchorcolor=blue, citecolor=blue,filecolor=blue,linkcolor=blue,menucolor=blue]{hyperref}

\begin{document}
\title{Random Recursive Hypergraphs}

\author{P. L. Krapivsky}
\affiliation{Department of Physics, Boston University, Boston, Massachusetts 02215, USA}
\affiliation{Santa Fe Institute, Santa Fe, New Mexico 87501, USA}

\begin{abstract} 
Random recursive hypergraphs grow by adding, at each step, a vertex and an edge formed by joining the new vertex to a randomly chosen existing edge. The model is parameter-free, and several characteristics of emerging hypergraphs admit neat expressions via harmonic numbers, Bernoulli numbers, Eulerian numbers, and Stirling numbers of the first kind. Natural deformations of random recursive hypergraphs give rise to fascinating models of growing random hypergraphs. 
\end{abstract}

\maketitle

\section{Introduction}

A hypergraph is a pair $(V,E)$ where $V$ is a set of vertices and  $E$ is a set of edges. An edge in a hypergraph is a non-empty subset of $V$. Thus the total number of edges is at most $2^{|V|}-1$. More precisely, this is valid for {\em simple} hypergraphs, i.e., hypergraphs without repeated edges; we will consider only simple hypergraphs. 

Hypergraphs \cite{Berge,HG-book,Jost22} provide a natural extension of graphs \cite{Diestel,Flajolet} and simplicial complexes \cite{Hatcher,Dima15,Ginestra16}. Hypergraphs and similar objects known as multilayer and higher-order networks \cite{Ginestra18,Bick21,Soumen22} encode higher-order interactions \cite{Shang22}. Such interactions are inevitable in several ecological, financial, transportation and social networks \cite{Harary56,Davis67,KR03,KR05-balance,KR06-balance,Kleinberg09,Kleinberg10,Pretolani,Leskovec16,Ghrist16,Lepri21,Ghrist21,Petri22a,Kleinberg23,Kleinberg22b}, and play important role in brain networks \cite{brain05,brain09,brain11}. Binary interactions are traditionally used in physics, but higher-order interactions appear, e.g., in recent studies of toy models of quantum chaos \cite{chaos20,chaos21a,chaos21b,chaos21c,chaos22}. 

Random graphs are well-explored \cite{Drmota,Newman-book,Frieze,Hofstad}. Studies of random simplicial complexes (see \cite{Pippenger,Linial06,Meshulam,Linial16,Farber16,Pittel16,Ginestra16-NGF,Farber17a,Farber17b,Bobrowski17,Ginestra17,Ginestra-Sergey,Ginestra18-NGF,Kahle,Petri-B21,Petri-B22,Dima22} and references therein) are more recent. Simplicial complexes are beautiful but they require too strict mutual inclusion of interactions. Hypergraphs relax the assumption of mutual inclusion and represent a much broader class of systems with higher-order interactions. The analyses of random hypergraphs are gaining popularity. As with other large random structures, static random hypergraphs have been the first research subject \cite{Shamir,Newman09}, and they are still actively investigated \cite{Chodrow20,Dumitriu21b,Masuda,Petri22b,Ren22,Marc22}. Special families of static random hypergraphs, e.g., regular (where each vertex belongs to the same number of edges) and uniform (where each edge contains the same number of vertices) are relatively well understood, see \cite{Shamir,Cooley18,Cooley21,Dumitriu21a,Li21,Isaev}. Statistical physics models on random hypergraphs is another growing area of research, see e.g.  \cite{Lenka17,Louise19,Louise20,Ginestra21,Cooley22}. 

Many hypergraphs are evolving, and our goal is to investigate a `null' model of growing random hypergraphs. We will ignore degradation, i.e., disappearance of vertices and edges, and consider hypergraphs growing via stochastic rules. When the hypergraphs become large, basic characteristics of these random hypergraphs are usually self-averaging, so the average values provide the chief information. (Some more subtle characteristics may remain non-self-averaging.) The cumulants and full probability distributions remain interesting even when a random quantity exhibits a self-averaging behavior, and for a few self-averaging quantities we will compute all cumulants and probability distributions. 

To motivate our model we begin with primordial hypergraph consisting of a single vertex and a single edge composed of that vertex. We study hypergraphs growing from the primordial hypergraph according to the following recursive rule: At each step, a new node $v$ is added together with a new edge obtained by adding $v$ to an edge $e$ chosen uniformly among existing edges. Hence the new edge is $e\cup {v}$. The number of vertices is always equal to the number of edges:
\begin{equation}
|V| = |E| = N
\end{equation}
This recursive procedure generates random recursive hypergraphs (RRHs). We call $N$ the size of the RRH. Treating $N$ as a (discrete) time variable allows one to study the evolution of RRHs. 

For $N=2$, the RRH hypergraph has the edge set
\begin{eqnarray*}
E = [\{v_1\}, ~\{v_1, v_2\}]
\end{eqnarray*}
Thus, the outcome of the growth procedure is deterministic for $N=2$. Starting from $N=3$, more than one hypergraphs can be built by the RRH procedure. For $N=3$, two edge sets 
\begin{subequations}
\begin{align}
\label{122}
&E = [\{v_1\}, ~\{v_1, v_2\}\,, ~\{v_1, v_3\}] \\
\label{123}
&E = [\{v_1\}, ~\{v_1, v_2\}\,, ~\{v_1, v_2,v_3\}]
\end{align}
\end{subequations}
are formed with equal probability. 

The RRHs model is parameter-free and its definition mimics the definition of random recursive trees (RRTs), a paradigmatic parameter-free model of growing random trees \cite{Drmota,Frieze,Hofstad}.  RRTs are engaging in their own right (see \cite{Pittel94,KR02,Janson05,Janson15,Janson19} and references therein). More importantly, the RRT is a mother model of growing networks as its simple deformations  lead to interesting growing network models \cite{Kleinberg,KR01,KR02-fluct,Slava05a,KR05,BK10,KR-copy,KR14,Bertoin15,KR17,Biham19a, Biham19b,Sumpter22}. For instance, adding redirection generates preferential attachment in some models \cite{Kleinberg,KR01,KR14} and non-self-averaging behaviors  in others \cite{KR14,KR17,Sumpter22}. Other simple deformations of the RRT generate dense networks \cite{KR-copy} that are poorly understood. 

Similarly to RRTs, hypergraphs built via the RRH rule are tractable. We now give a small sample of exact results derived in this paper. For instance, the total number $\mathcal{N}_1(N)$ of vertexes each belonging to a singe edge  admits a neat analytical description. We have $\mathcal{N}_1(3)=2$ for the hypergraph \eqref{122} and $\mathcal{N}_1(3)=1$ for the hypergraph \eqref{123}, and generally $\mathcal{N}_1(N)$ has the probability distribution
$P_1(n,N) = \text{Prob}[\mathcal{N}_1(N)=n]$ expressible via Eulerian numbers \cite{Euler36,Euler55,Knuth}:
\begin{equation}
\label{N1:Euler}
P_1(n,N)=\frac{1}{(N-1)!}\,\eulerian{N-1}{n-1}
\end{equation}

Another interesting random quantity is the total number $\mathcal{R}_2(N)$ of edges of size two. We have $\mathcal{R}_2(3)=2$ for the hypergraph \eqref{122} and $\mathcal{R}_2(3)=1$ for the hypergraph \eqref{123}. We will show that the random variable $\mathcal{R}_2(N)$ has the probability distribution $\Pi(r,N)=\text{Prob}[\mathcal{R}_2(N)=r]$ expressible via Stirling numbers \cite{Stirling,Knuth} of the first kind:
\begin{equation}
\label{R2:Stirling}
\Pi(r,N)=\frac{1}{(N-1)!}\,{N-1\brack r}
\end{equation}

In Sec.~\ref{sec:DD}, we study the degree distribution of RRHs. First, we derive the probability distribution  \eqref{N1:Euler} for the number of vertexes of degree one. We then demonstrate that the average degree $\langle \mathcal{N}_k(N)\rangle$ distribution exhibits remarkable stationarity, namely, it is strictly linear in $N$ for all degrees $k\leq N-1$. The edge size distribution of RRHs is studied in Sec.~\ref{sec:rank} where among other results we derive \eqref{R2:Stirling}. In Sec.~\ref{sec:leader}, we discuss leaders in degree \footnote{The primordial vertex has the highest degree $N$, so the leader is the vertex with the highest degree among vertexes different from the primordial.}. We also analyze the dependence of the degree distribution of a vertex on its index. For instance, the degree of the second vertex has a uniform degree distribution, so it is a non-self-averaging quantity strongly fluctuating from realization to realization. In Sec.~\ref{sec:leaves}, we study leaves and determine their probability distribution in the $N\to\infty$ limit. In Sec.~\ref{sec:RD}, we deform the RRH model.  Specifically, we employ the same growth rule as in the RRH model with probability $1-r$; with probability $r$, we redirect the randomly chosen edge $e$ to its maternal edge $e'$ and call $e'\cup v$ the new edge. In Sec.~\ref{sec:concl}, we comment on possible directions of future research and briefly mention a few natural deformations of the RRH model.

\section{Degree Distribution}
\label{sec:DD}

The degree of a vertex in a hypergraph is the number of edges containing it
\begin{equation}
\label{degree:def}
d(v) = \#\{e\,| ~v\in e\}
\end{equation}
Let $\mathcal{N}_k(N)$ be the total number of vertices of degree $k$:
\begin{equation}
\label{Nk:def}
\mathcal{N}_k= \#\{v: \,d(v)=k\}
\end{equation}
The RRHs of size $N=1$ and $N=2$ are deterministic: $\mathcal{N}_1(1)=1$ and $\mathcal{N}_1(2)=\mathcal{N}_2(2)=1$. For $N \geq 3$, the quantities $\mathcal{N}_k(N)$ with $k<N$ are random. Generally 
\begin{equation}
\label{norm}
\sum_{k=1}^N \mathcal{N}_k(N) = N
\end{equation}
The maximal degree is equal to the size of the hypergraph, and there is only one vertex with maximal degree, the primordial vertex that belongs to every edge:
\begin{equation}
\label{max-degree}
\mathcal{N}_N(N) = 1
\end{equation}

\subsection{Vertices of degree one}
\label{subsec:1}

The total number of vertices of the smallest degree, $\mathcal{N}_1(N)$, is a random quantity for $N\geq 3$. We already know
\begin{equation}
\mathcal{N}_1(3) = 
\begin{cases}
2 & \text{prob} ~~1/2\\
1 & \text{prob} ~~1/2
\end{cases}
\end{equation}
The recursive nature of the process shows that 
\begin{equation}
\label{N1N}
\mathcal{N}_1(N+1) = 
\begin{cases}
\mathcal{N}_1+1 & \text{prob} ~~1-\mathcal{N}_1/N\\
\mathcal{N}_1    & \text{prob} ~~ \mathcal{N}_1/N
\end{cases}
\end{equation}
for $N\geq 2$. Hereinafter we often omit the dependence of $N$, e.g., $\mathcal{N}_1$ means $\mathcal{N}_1(N)$; we explicitly write $\mathcal{N}_1(N+1)$ to avoid confusion. Averaging \eqref{N1N} shows that $N_1=\langle \mathcal{N}_1\rangle$ satisfies the recurrence
\begin{equation}
\label{N1-rec}
N_1(N+1)=\left(1-\frac{1}{N}\right)N_1+1
\end{equation}
for $N\geq 2$. Using the known value $N_1(2)=1$ as the boundary condition we solve \eqref{N1-rec} and find
\begin{equation}
\label{N1-sol}
\langle \mathcal{N}_1\rangle = N_1 = \frac{N}{2}\,, \qquad N\geq 2
\end{equation}

Similarly one can compute higher moments of the random quantity $\mathcal{N}_1$. The second moment $\langle \mathcal{N}_1^2\rangle$ satisfies the recurrence 
\begin{equation}
\label{N2-rec}
\langle \mathcal{N}_1^2(N+1)\rangle=\left(1-\frac{2}{N}\right)\langle \mathcal{N}_1^2\rangle+N+\frac{1}{2}
\end{equation}
when $N\geq 2$. Solving \eqref{N2-rec} subject to the boundary condition $\langle \mathcal{N}_1^2(3)\rangle=5/2$ we find that for $N\geq 3$
\begin{equation}
\label{N2-sol}
\langle \mathcal{N}_1^2\rangle= \frac{N(3N+1)}{12}
\end{equation}
The variance $V_1=\langle \mathcal{N}_1^2\rangle_c\equiv \langle [\mathcal{N}_1-N_1]^2\rangle$ is 
\begin{equation}
\label{Var-N1}
V_1= \frac{N}{12}
\end{equation}
for $N\geq 3$. 

The random quantity $\mathcal{N}_1(N)$ is concentrated around its average in the large $N$ limit. More precisely, the probability distribution 
$P_1(n,N)=\text{Prob}[\mathcal{N}_1(N)=n]$ is asymptotically Gaussian:
\begin{equation}
\label{P1nN:Gauss}
P_1(n, N)\simeq \sqrt{\frac{6}{\pi N}}\,\text{exp}\!\left[-\frac{6(n -N/2)^2}{N}\right]
\end{equation}
This assertion can be deduced from the expression the probability distribution $P_1(n,N)$ via Eulerian numbers [Eq.~\eqref{N1:Euler} which we derive below]. The limit law \eqref{P1nN:Gauss} then follows from \eqref{N1:Euler}, see \cite{Eulerian} for details. The convergence rate to the limit law \eqref{P1nN:Gauss} is also known \cite{Eulerian}. 

The Gaussian limiting behavior \eqref{P1nN:Gauss} becomes intuitively plausible after realizing that higher cumulants grow anomalously slow with $N$ (odd cumulants even vanish). To appreciate this assertion let us look first at the third moment $\langle \mathcal{N}_1^3\rangle$. It obeys the recurrence
\begin{equation*}
\label{N3-rec}
\langle \mathcal{N}_1^3(N+1)\rangle=\left(1-\frac{3}{N}\right)\langle \mathcal{N}_1^3\rangle+\frac{(N+1)(3N+1)}{4}
\end{equation*}
which is solved for $N\geq 4$ to find
\begin{equation}
\label{N3-sol}
\langle \mathcal{N}_1^3\rangle= \frac{N^2(N+1)}{8}
\end{equation}
The prediction \eqref{N3-sol} remains valid even for $N=3$. Surprisingly, the third cumulant vanishes
\begin{eqnarray}
\label{N1-cum3}
\langle \mathcal{N}_1^3\rangle_c\equiv \langle [\mathcal{N}_1-N_1]^3\rangle = 0
\end{eqnarray}
The fourth moment $\langle \mathcal{N}_1^4\rangle$ obeys the recurrence
\begin{equation*}
\label{N4-rec}
\langle \mathcal{N}_1^4(N+1)\rangle=\left(1-\frac{4}{N}\right)\langle \mathcal{N}_1^4\rangle+\frac{6N^3+15N^2+9N+2}{12}
\end{equation*}
which is solved to give 
\begin{equation}
\label{N4-sol}
\langle \mathcal{N}_1^4\rangle= \frac{N(15N^3+30N^2+5N-2)}{240}
\end{equation}
for $N\geq 5$. Combining \eqref{N4-sol} with previous results \eqref{N1-sol}, \eqref{N2-sol}, \eqref{N3-sol} for lower moments we extract a neat expression for the fourth cumulant:
\begin{equation}
\label{N1-cum4}
\langle \mathcal{N}_1^4\rangle_c = -\frac{N}{120}
\end{equation}
for $N\geq 5$. 

Equations \eqref{N1-sol}, \eqref{Var-N1}, \eqref{N1-cum3} and \eqref{N1-cum4} for the cumulants $\langle \mathcal{N}_1^p\rangle_c$ with $p=1,2,3,4$, suggest that all cumulants are strictly linear in $N$. More precisely, 
\begin{equation}
\label{N1-cump}
\langle \mathcal{N}_1^p\rangle_c = \kappa_p N
\end{equation}
is expected to hold for $N > p$. In other words, for any $N>p$, the fractions $\kappa_p=N^{-1}\langle \mathcal{N}_1^p\rangle_c$ are {\em stationary}, that is, independent on $N$. This observation is borne out of straightforward calculations for small $p$ which are difficult to extend to large $p$. 

We now derive the announced result \eqref{N1-cump}. Using \eqref{N1N} we deduce a recurrence
\begin{eqnarray}
\label{FnN:eq}
P_1(n+1,N+1) &=& \left(1-\frac{n}{N}\right)P_1(n,N)\nonumber \\
&+& \frac{n+1}{N}\,P_1(n+1,N)
\end{eqnarray}
Making the substitution 
\begin{equation}
\label{FnN:ansatz}
P_1(n+1,N+1)=\frac{1}{N!}\,\eulerian{N}{n}
\end{equation}
we recast \eqref{FnN:eq} into 
\begin{equation}
\label{Phi-nN:eq}
\eulerian{N}{n}=(N-n)\eulerian{N-1}{n-1}+(n+1)\eulerian{N-1}{n}
\end{equation}
and recognize that this neat recurrence is an addition formula for Eulerian numbers \cite{Euler36,Euler55,Knuth}. One can check that the boundary conditions agree with the standard definition of Eulerian numbers. 

Equation \eqref{FnN:ansatz} gives the probability distribution for the number of vertices of degree one.  Using explicit expressions \cite{Knuth} for Eulerian numbers 
\begin{equation*}
\begin{split}
\eulerian{N}{0} &=1 \\
\eulerian{N}{1}  &=2^N-N-1 \\
\eulerian{N}{2} &=3^N-(N+1)2^N+\binom{N+1}{2} \\
\eulerian{N}{3} &=4^N-(N+1)3^N+2^N\binom{N+1}{2}-\binom{N+1}{3}
\end{split}
\end{equation*}
we obtain $F(n,N+1)$ for $n=1,2,3,4$. The prediction $F(1,N+1)=1/N!$ is obvious, while a direct straightforward derivation of $F(n,N+1)$ for $n=2,3,4$ is laborious. 

Using the basic identity \cite{Knuth}
\begin{equation}
\label{mirror:E}
\eulerian{N}{n}=\eulerian{N}{N-1-n}
\end{equation}
reflecting the mirror symmetry between Eulerian numbers one obtains 
\begin{equation}
\label{mirror}
P_1(n+1,N+1) = P_1(N-n,N+1)
\end{equation}
giving $P_1(n,N+1)$ for $n=N, N-1,N-2,N-3$. 

Thanks to numerous identities \cite{Knuth,Kyle} satisfied by Eulerian numbers one can \eqref{N1-cump} and establish the amplitudes $\kappa_p$. This has been done in Refs.~\cite{David62,Janson13} where Eq.~\eqref{N1-cump} was derived and the amplitudes were expressed through Bernoulli numbers $B_p$:
\begin{equation}
\label{N1-cum-B}
\langle \mathcal{N}_1^p\rangle_c =p^{-1} B_p N
\end{equation}
This holds when $p<N$, see \cite{David62,Janson13}. Bernoulli numbers are the coefficients in the power series 
\begin{equation}
\label{B:def}
\frac{z}{e^z-1}+z =\sum_{p\geq 0}B_p\,\frac{z^p}{p!}
\end{equation}
Adding or not the second term on the left-hand side in \eqref{B:def} changes only $B_1$, and the choice of the best convention is the matter of debate. The definition \eqref{B:def} gives $B_1=\frac{1}{2}$, so Eq.~\eqref{N1-cum-B} remains valid when $p=1$. Furthermore, two useful formulas 
\begin{equation*}
\begin{split}
& \sum_{m=0}^n (-1)^m\, \eulerian{n}{m} = 2^{n+1}\big(2^{n+1}-1\big)\,\frac{B_{n+1}}{n+1}\\
& \sum_{m=0}^n (-1)^m\, \eulerian{n}{m} \binom{n}{m}^{-1} = (n+1)B_n
\end{split}
\end{equation*}
connecting Eulerian numbers to Bernoulli numbers are valid for $n\geq 0$ if $B_1=\frac{1}{2}$, while if $B_1$ is set to $-\frac{1}{2}$ they are applicable only when  $n\geq 1$ and $n\geq 2$, respectively.

Since $B_p=0$ when $p\geq 3$ is odd,  see \cite{Knuth}, we have 
\begin{equation*}
\langle \mathcal{N}_1^3\rangle_c = \langle \mathcal{N}_1^5\rangle_c =\ldots= \langle \mathcal{N}_1^{2\lfloor N/2\rfloor -1}\rangle_c=0
\end{equation*}
The first of this relations is \eqref{N1-cum3} which we derived above using straightforward calculations. Even cumulants are non-zero, and their signs alternate. Even cumulants following \eqref{Var-N1} and \eqref{N1-cum4} are
\begin{equation*}
\begin{split}
\langle \mathcal{N}_1^6\rangle_c &= \tfrac{N}{252}\,, \quad~\langle \mathcal{N}_1^8\rangle_c = - \tfrac{N}{240}\,, 
~~\langle \mathcal{N}_1^{10}\rangle_c = \tfrac{N}{132}\\
\langle \mathcal{N}_1^{12}\rangle_c &= -\tfrac{691 N}{32760}\,, ~~\langle \mathcal{N}_1^{14}\rangle_c = \tfrac{N}{12}\,, 
~~\langle \mathcal{N}_1^{16}\rangle_c = \tfrac{3617N}{8160}
\end{split}
\end{equation*}
etc.

\subsection{Vertices of higher degree: Average degree distribution}
\label{subsec:k}

When $k\geq 2$, the random quantity $\mathcal{N}_k$ evolves according to stochastic rule 
\begin{equation}
\label{NkN}
\mathcal{N}_k(N+1) =
\begin{cases}
\mathcal{N}_k+1  &\text{prob} ~~ \frac{(k-1)\mathcal{N}_{k-1}}{N}\\
\mathcal{N}_k-1  &\text{prob} ~~ \frac{k\mathcal{N}_{k}}{N} \\
\mathcal{N}_k     &\text{prob} ~~  1-\frac{(k-1)\mathcal{N}_{k-1}+k \mathcal{N}_k}{N}
\end{cases}
\end{equation}

Averaging \eqref{NkN} we find that $N_k=\langle \mathcal{N}_k\rangle$ satisfies 
\begin{equation}
\label{Nk-rec}
N_k(N+1)=\left(1-\frac{k}{N}\right)N_k + \frac{k-1}{N}\,N_{k-1}
\end{equation}
This recurrence admits a remarkably simple solution
\begin{equation}
\label{Nk-sol}
N_k= \frac{N}{k(k+1)}\,, \qquad k \leq N-1
\end{equation}
Note that  
\begin{equation}
\label{check}
\sum_{k=1}^N N_k = \sum_{k=1}^{N-1}\frac{N}{k(k+1)}+1=N
\end{equation}
where we have used \eqref{Nk-sol} and \eqref{max-degree}. Equation \eqref{check} should be valid due to the exact sum rule \eqref{norm}, so confirming it provides a consistency check. 

We emphasize that there are no sub-leading terms in \eqref{Nk-sol}, that is, the average number of vertices $N_k(N)$ of degree $k$ exhibits a strictly linear in $N$ behavior. Thus for any $N>k$, the fractions $n_k=N_k(N)/N$ are {\em stationary}, that is, independent on $N$:
\begin{equation}
\label{nk-sol}
n_k= \frac{1}{k(k+1)}
\end{equation}

For RRTs, the asymptotic behavior of the degree distribution is exponential, $n_k = 2^{-k}$, but sub-leading terms do not vanish \cite{KR02-fluct}. Interestingly, for RRHs, one gets algebraic behavior \eqref{nk-sol} without preferential attachment. The lack of correction terms in the degree distribution of the RRH is striking:  For all degrees, $k\leq N-1$, the behavior of the degree distribution of the RRH is such as if the system was effectively infinite. 

\subsection{Vertices of degree two: Fluctuations}
\label{subsec:2-fluct}

The computations of fluctuations of the random quantities $\mathcal{N}_k$ become more involved as $k$ increases. Here we consider $\mathcal{N}_2$. Specializing \eqref{NkN} to $k=2$ gives
\begin{equation}
\label{N2N}
\mathcal{N}_2(N+1) =
\begin{cases}
\mathcal{N}_2+1  &\text{prob} ~~ \frac{\mathcal{N}_1}{N}\\
\mathcal{N}_2-1  &\text{prob} ~~ \frac{2\mathcal{N}_2}{N} \\
\mathcal{N}_2     &\text{prob} ~~  1-\frac{\mathcal{N}_1+2\mathcal{N}_2}{N}
\end{cases}
\end{equation}
Let us try to determine the variance. Taking the square of \eqref{N2N} and averaging we derive the governing equation for the second moment  $\langle \mathcal{N}_2^2\rangle$: 
\begin{equation}
\label{N22}
\langle \mathcal{N}_2^2(N+1)\rangle=\left(1-\frac{4}{N}\right)\langle \mathcal{N}_2^2\rangle+\frac{2}{N}\,\langle \mathcal{N}_1\mathcal{N}_2\rangle+\frac{5}{6}
\end{equation}
Equation \eqref{N22} shows that we need an additional equation for the correlation function $\langle \mathcal{N}_1\mathcal{N}_2\rangle$. We first note that the product $\mathcal{N}_1(N+1)\mathcal{N}_2(N+1)$ evolves according to 
\begin{equation}
\label{N12N}
\begin{cases}
\mathcal{N}_1(\mathcal{N}_2+1)  &\text{prob} ~~ \frac{\mathcal{N}_{1}}{N}\\
(\mathcal{N}_1+1)(\mathcal{N}_2-1)  &\text{prob} ~~ \frac{2\mathcal{N}_{2}}{N} \\
(\mathcal{N}_1+1)\mathcal{N}_2     &\text{prob} ~~  1-\frac{\mathcal{N}_{1}+2\mathcal{N}_{2}}{N}
\end{cases}
\end{equation}
Averaging  \eqref{N12N} and taking into account \eqref{N2-sol} and \eqref{Nk-sol} with $k=1, 2$ we obtain
\begin{eqnarray}
\label{N12}
\langle \mathcal{N}_1(N+1)\mathcal{N}_2(N+1)\rangle &=& \left(1-\frac{3}{N}\right)\langle \mathcal{N}_1\mathcal{N}_2\rangle \nonumber \\
&+&  \frac{5N-3}{12}
\end{eqnarray}
A straightforward calculation gives
\begin{equation}
\label{N123}
(\mathcal{N}_1, \mathcal{N}_2, \mathcal{N}_3) =
\begin{cases}
(3,0,0)  &\text{prob} ~~ \frac{1}{6}  \\
(2,1,0)  &\text{prob} ~~ \frac{1}{2} \\
(2,0,1)  &\text{prob} ~~ \frac{1}{6} \\
(1,1,1)  &\text{prob} ~~ \frac{1}{6} 
\end{cases}
\end{equation}
for $N=4$. Thus $\langle \mathcal{N}_1\mathcal{N}_2\rangle=\frac{7}{6}$ when $N=4$. Using this value as a boundary condition we solve the 
 recurrence \eqref{N12} to yield
\begin{equation}
\label{N12-sol}
\langle \mathcal{N}_1\mathcal{N}_2\rangle=\frac{N(N-1)}{12} +\frac{1}{(N-1)(N-2)(N-3)}
\end{equation}
for $N\geq 4$. Using \eqref{N12-sol} we find that the centered pair correlation $\langle \mathcal{N}_1\mathcal{N}_2\rangle_c=\langle \mathcal{N}_1\mathcal{N}_2\rangle-\langle \mathcal{N}_1\rangle \langle\mathcal{N}_2\rangle$ is given by
\begin{equation}
\label{N12-c}
\langle \mathcal{N}_1\mathcal{N}_2\rangle_c=-\frac{N}{12}+\frac{1}{(N-1)(N-2)(N-3)}
\end{equation}

Inserting \eqref{N12-sol} into \eqref{N22} gives a closed recurrence for $\langle \mathcal{N}_2^2\rangle$. Solving this recurrence yields
\begin{equation}
\label{N22-sol}
\langle \mathcal{N}_2^2\rangle=\frac{N(5N+23)}{180} + \frac{1}{(N-1)(N-2)(N-3)}
\end{equation}
for $N\geq 5$. The variance $V_2=\langle [\mathcal{N}_2-N_2]^2\rangle$  reads
\begin{equation}
\label{Var-N2}
V_2=\frac{23N}{180} + \frac{1}{(N-1)(N-2)(N-3)}
\end{equation}

The random quantity $\mathcal{N}_2(N)$ apparently concentrates around its average when $N\to\infty$. More precisely, the probability distribution $P_2(n,N)=\text{Prob}[\mathcal{N}_2(N)=n]$ is believed to be asymptotically Gaussian:
\begin{equation}
P_2(n,N) \simeq \sqrt{\frac{90}{23\pi N}}\,\text{exp}\!\left[-\frac{90(n -N/6)^2}{23 N}\right]
\end{equation}

Generally $\langle \mathcal{N}_i\mathcal{N}_j\rangle_c=\langle \mathcal{N}_i\mathcal{N}_j\rangle-\langle \mathcal{N}_i\rangle \langle\mathcal{N}_j\rangle$ for all $i\leq j$ are expected to grow linearly with $N$:
\begin{equation}
\lim_{N\to\infty}N^{-1}\langle \mathcal{N}_i\mathcal{N}_j\rangle_c=\nu_{i,j}
\end{equation}

The quantities $\nu_{i,j}$ appear to be rational numbers, and $\nu_{i,i}>0$ for all $i\geq 1$ as the variances are positive. The above {\em exact} calculations [viz., Eqs.~\eqref{Var-N1}, \eqref{N12-c}, \eqref{Var-N2}] give
\begin{equation}
\nu_{1,1}=\frac{1}{12}, \quad \nu_{1,2}=-\frac{1}{12}, \quad \nu_{2,2}=\frac{23}{180}
\end{equation}

Finding $\nu_{i,j}$ for all $i\leq j$ is a challenge. For RRTs, similar calculations have been performed (see, e.g., \cite{Janson05}), but no simple formulas giving $\nu_{i,j}$ for all $i\leq j$ have been found to the best of our knowledge. 

\section{Rank Distribution}
\label{sec:rank}

The rank of a vertex in a hypergraph is the size of the minimal edge containing it
\begin{equation}
r(v) = \text{min}\{|e|: v\in e\}
\end{equation}
For instance, for a hypergraph with edge set
\begin{subequations}
\label{14}
\begin{equation}
\label{E14}
\begin{split}
 & \{v_1\} \\
 & \{v_1, v_2\}, ~ \{v_1, v_3\}, ~\{v_1, v_8\}, ~ \{v_1, v_9\},~\{v_1, v_{10}\} \\
 & \{v_1, v_2, v_4\}, ~\{v_1, v_3,v_5\}, ~\{v_1, v_{10},v_{12}\}\\
 & \{v_1, v_3,v_5,v_6\}, ~\{v_1, v_2, v_4,v_7\}, ~  \{v_1, v_2, v_4, v_{14}\}\\
 & \{v_1, v_3,v_5,v_6,v_{11}\}, ~\{v_1, v_2, v_4,v_7,v_{13}\}
 \end{split}
\end{equation}
the ranks of the vertices are
\begin{equation}
\label{r14}
\begin{split}
 & r(v_1)=1 \\
 & r(v_2) = r(v_3) = r(v_8) = r(v_9) = r(v_{10})=2  \\
 & r(v_4) = r(v_5) = r(v_{12}) = 3   \\
 & r(v_6) = r(v_7) = r(v_{14}) = 4   \\
 & r(v_{11}) = r(v_{13}) = 5
 \end{split}
\end{equation}
\end{subequations}

Let $\mathcal{R}_k$ be the total number of vertices of rank $k$:
\begin{equation}
\mathcal{R}_k= \#\{v: \,r(v)=k\}
\end{equation}
For the hypergraph  \eqref{14},
\begin{equation*}
\mathcal{R}_1=1, \quad \mathcal{R}_2=5, \quad \mathcal{R}_3 = \mathcal{R}_4=3, \quad \mathcal{R}_5=2
\end{equation*}
Generally for any hypergraph of size $N$
\begin{subequations}
\begin{equation}
\label{norm-R}
\sum_{k=1}^N \mathcal{R}_k = N
\end{equation}
For RRHs, there is exactly one vertex with minimal rank, the primordial vertex, see \eqref{r14}. Thus
\begin{equation}
\label{R1}
\mathcal{R}_1 = 1
\end{equation}
The maximal possible rank of a hypergraph is equal to its size. In contrast to the minimal rank which is always realized by the definition of RRHs, only the last vertex $v_N$ in an RRH may have the maximal rank $N$. The recursive building procedure of the RRHs implies that the vertex with rank $N$ arisses with probability $1/(N-1)!$. Hence the average number of vertices of the maximal possible rank, $R_N=\langle \mathcal{R}_N\rangle$, is
\begin{equation}
\label{RN}
R_N = \frac{1}{(N-1)!}
\end{equation}
\end{subequations}

Note that the number of vertices of rank $k$ is equal to the number of edges $\mathcal{E}_k$ of size $k$:
\begin{equation}
\mathcal{E}_k = \#\{e: \, |e|=k\} = \mathcal{R}_k
\end{equation}

The recursive nature of the RRHs leads to the following stochastic evolution equation for $\mathcal{R}_k$:
\begin{equation}
\label{RkN}
\mathcal{R}_k(N+1) = 
\begin{cases}
\mathcal{R}_k+1 & \text{prob} ~~\mathcal{R}_{k-1}/N\\
\mathcal{R}_k    & \text{prob} ~~ 1-\mathcal{R}_{k-1}/N
\end{cases}
\end{equation}
Indeed, random quantities $\mathcal{R}_j$ can only increase, and adding a vertex increases one of $\mathcal{R}_j$ by one [cf. Eq.~\eqref{norm-R}]. The probability that a new vertex joins an edge of size $k-1$ is $\mathcal{E}_{k-1}/N=\mathcal{R}_{k-1}/N$, and if this happens the new vertex has rank $k$. The recurrence \eqref{RkN} is valid for $k=2,\ldots,N$, and even for $k=1$ and $k=N+1$ if we recall that $\mathcal{R}_0=0$ and $\mathcal{R}_{N+1}(N)=0$.  

\subsection{Vertices of rank two}
\label{subsec:2}

When $k=2$, the recurrence \eqref{RkN} becomes
\begin{equation}
\label{R2N}
\mathcal{R}_2(N+1) = 
\begin{cases}
\mathcal{R}_2+1 & \text{prob} ~~1/N\\
\mathcal{R}_2    & \text{prob} ~~ 1-1/N
\end{cases}
\end{equation}
Averaging \eqref{R2N} shows that $R_2=\langle \mathcal{R}_2\rangle$ satisfies the recurrence
\begin{equation}
\label{R2-rec}
R_2(N+1) = R_2 + \frac{1}{N}
\end{equation}
from which
\begin{equation}
\label{R2-sol}
R_2(N) = H_{N-1}
\end{equation}
where $H_n=\sum_{1\leq i\leq n} i^{-1}$ are harmonic numbers \cite{Knuth}. Using  \eqref{R2N} one similarly deduces the recurrence for the second moment 
\begin{equation}
\label{R22-rec}
\langle \mathcal{R}_2^2(N+1)\rangle = \langle \mathcal{R}_2^2\rangle  +  \frac{2}{N}\,H_{N-1}+\frac{1}{N}
\end{equation}
from which
\begin{equation}
\label{R22-sol}
\langle \mathcal{R}_2^2\rangle  = (H_{N-1})^2 + H_{N-1}-H_{N-1}^{(2)}
\end{equation}
This equation and exact results for higher moments $\langle \mathcal{R}_2^p\rangle$ involve generalized harmonic numbers:
\begin{equation}
\label{Harmonic}
H_n^{(p)}=\sum_{i=1}^n \frac{1}{i^p}
\end{equation}
Combining \eqref{R2-sol} and \eqref{R22-sol} we find the variance of the number of vertices of rank two:
\begin{equation}
\label{R2-var}
\langle \mathcal{R}_2^2\rangle - \langle \mathcal{R}_2\rangle^2  = H_{N-1}-H_{N-1}^{(2)}
\end{equation}

We now show how to compute the entire probability distribution $\Pi(r,N)=\text{Prob}[\mathcal{R}_2(N)=r]$ for the number of vertices of rank two.  Using \eqref{R2N} we deduce
\begin{equation}
\label{PrN:eq}
\Pi(r,N+1)=\frac{1}{N}\,\Pi(r-1,N)+\left(1-\frac{1}{N}\right)\Pi(r,N)
\end{equation}
Making the substitution \eqref{R2:Stirling} we recast \eqref{PrN:eq} into
\begin{equation}
\label{PirN:eq}
{N\brack r}={N-1\brack r-1}+(N-1){N-1\brack r}
\end{equation}
which is an addition formula for Stirling numbers of the first kind \cite{Knuth}. This completes the derivation of the announced formula \eqref{R2:Stirling}. 

Using well-known expressions \cite{Knuth} for the extremal Stirling numbers of the first kind 
\begin{subequations}
\begin{equation}
\label{S:extremal}
{n\brack 1} = (n-1)!, \qquad {n\brack n} = 1
\end{equation}
one gets
\begin{equation}
\label{PrN:extremal}
P(1,N)=\frac{1}{N-1}\,, \qquad P(N-1,N)=\frac{1}{(N-1)!}
\end{equation}
\end{subequations}
These values readily follow directly from the definition of the RRH. Using expressions \cite{Knuth} for penultimate extremal Stirling numbers of the first kind 
\begin{subequations}
\begin{equation}
\label{S:extremal-2}
{n\brack 2} = (n-1)! H_{n-1}, \qquad {n\brack n-1} = \binom{n}{2}
\end{equation}
one deduces
\begin{equation}
\label{PrN:extremal-2}
P(2,N)=\frac{H_{N-2}}{N-1}\,, \quad P(N-2,N)=\frac{1}{2(N-3)!}
\end{equation}
\end{subequations}
which are harder to derive in a straightforward manner, namely merely relying on the definition of the RRH.

\subsection{Vertices of higher ranks}
\label{subsec:Rk}

Averaging \eqref{RkN} leads to the recurrence
\begin{equation}
\label{Rk-rec}
R_k(N+1) = R_k + \frac{1}{N}\,R_{k-1}
\end{equation}
Using the recursive nature of Eqs.~\eqref{Rk-rec} we begin with $R_1=1$ and solve \eqref{Rk-rec} for all $k$. The general result is
\begin{equation}
\label{Rk-sol}
R_{k+1}(N+1) = \sum_{1\leq j_1<\ldots<j_k\leq N}\frac{1}{j_1\times \ldots\times j_k}
\end{equation}
All restrictions on the sum are indicated in \eqref{Rk-sol}. When $k>N$, the sum is empty leading to $R_{k+1}(N+1)=0$. When $k=N$, the sum contains a single term with $j_i=i$ and hence  $R_{N+1}(N+1)=1/N!$ in agreement with \eqref{RN}. As another consistency check we note that \eqref{Rk-sol} at $k=1$ reduces  to \eqref{R2-sol}. Specializing \eqref{Rk-sol} to $k=2$ and expressing the sum via harmonic numbers yields 
\begin{equation}
\label{R3-sol}
R_3(N+1) = \frac{(H_{N})^2 -H_{N}^{(2)}}{2}
\end{equation}

Simplifying the sum in \eqref{Rk-sol} and expressing the exact solution via known finite sums like harmonic numbers is feasible. However, the results become more and more cumbersome as $k$ increases. Here we only mention the leading asymptotic behavior for $k=O(1)$ and $N\to\infty$. In this situation, one can replace the recurrence \eqref{Rk-rec} by the differential equation
\begin{equation}
\label{Rk-DE}
\frac{dR_k}{dN} = \frac{1}{N}\,R_{k-1}
\end{equation}
Solving these equations recurrently starting from $R_1=1$, see \eqref{R1}, yields
\begin{equation}
\label{Rk-asymp}
R_{k+1} = \frac{(\ln N)^k}{k!}
\end{equation}

Let us gauge the accuracy of Eq.~\eqref{Rk-asymp}. First,  we recall the asymptotic formulas \cite{Knuth}
\begin{subequations}
\begin{align}
\label{Harmonic-1}
H_n          & = \ln n + \gamma +\frac{1}{2n} + O(n^{-2}) \\
\label{Harmonic-p}
H_n^{(p)} & = \zeta(p)-\frac{1}{(p-1)n^{p-1}} + O(n^{-p}) 
\end{align}
\end{subequations}
Here $\gamma=0.57721\ldots$ is the Euler constant, $\zeta(\cdot)$ is the zeta function, and \eqref{Harmonic-p} is valid for integer $p\geq 2$. Using \eqref{R2-sol} and \eqref{Harmonic-1} we obtain
\begin{subequations}
\begin{equation}
\label{R2-asymp}
R_2 = \ln N + \gamma - \frac{1}{2N} + O(N^{-2})
\end{equation}
Using \eqref{R3-sol} and \eqref{Harmonic-p} with $p=2$ we obtain
\begin{eqnarray}
\label{R3-asymp}
R_3 &=& \tfrac{1}{2}(\ln N)^2 +\gamma \ln N + \tfrac{1}{2}\gamma^2 - \tfrac{1}{12}\pi^2 \nonumber\\
        &+& \frac{1-\gamma \ln N}{N}+O(N^{-2})
\end{eqnarray}
\end{subequations}
The leading terms in \eqref{R2-asymp}--\eqref{R3-asymp} agree with the general leading asymptotic predicted by Eq.~\eqref{Rk-asymp}. 

The asymptotic expansions \eqref{R2-asymp}--\eqref{R3-asymp} also suggest an improvement of the leading asymptotic \eqref{Rk-asymp}, e.g., next two  sub-leading terms in the expansion  
\begin{eqnarray}
\label{Rk-asymp-2}
R_{k+1} &=& \frac{(\ln N)^k}{k!} + \gamma\,\frac{(\ln N)^{k-1}}{(k-1)!}  \nonumber \\
&+& \frac{6\gamma^2 - \pi^2}{12}\,\frac{(\ln N)^{k-2}}{(k-2)!}+\ldots
\end{eqnarray}
are easy to confirm by analyzing \eqref{Rk-rec}.

\section{Leaders}
\label{sec:leader}

In RRTs and similar growing networks, the statistical properties of the node of the highest degree can be rather remarkable \cite{KR02,JML08}. In RRHs, the primordial vertex has the highest possible degree, so it is natural to define the leader in degree as the vertex with the highest degree among vertices different from the primordial. 

Take the second vertex, the most plausible candidate for having the second highest degree. The probability distribution $P_2(d,N)=\text{Prob}\{d[v_2(N)]=d\}$ of the degree of the second vertex satisfies the recurrence
\begin{equation}
\label{PdN-2}
P_2(d,N+1)=\tfrac{d-1}{N}P_2(d-1,N)+\tfrac{N-d}{N}P_2(d,N)
\end{equation}
for $N\geq 2$. Using $P_2(1,2)=1$ as initial condition and iterating \eqref{PdN-2} we arrive at a very simple solution
\begin{subequations}
\label{PdN:sols}
\begin{equation}
\label{PdN:sol-2}
P_2(d,N)=\frac{1}{N-1}\,, \qquad d=1,\ldots,N-1
\end{equation}
Thus the probability distribution $P_2(d,N)$ is uniform, so the random quantity $d=d[v_2(N)]$ is non-self-averaging. 

For the third vertex,  the probability distribution of its degree, $P_3(d,N)=\text{Prob}\{d[v_3(N)]=d\}$, satisfies the same equation as $P_2(d,N)$, viz.
\begin{equation*}
P_3(d,N+1)=\tfrac{d-1}{N}P_3(d-1,N)+\tfrac{N-d}{N}P_3(d,N)
\end{equation*}
which is valid for $N\geq 3$. The initial condition is different, $P_3(1,3)=1$, so the solution also differs from $P_2(d,N)$: 
\begin{equation}
\label{PdN:sol-3}
P_3(d,N)=\frac{2(N-d-1)}{(N-1)(N-2)}
\end{equation}
In contrast to \eqref{PdN:sol-2}, the distribution \eqref{PdN:sol-3} is not uniform.

Generally for the $m^\text{th}$ vertex, the  probability distribution $P_m(d,N)=\text{Prob}\{d[v_m(N)]=d\}$ satisfies an equation mathematically identical to \eqref{PdN-2}. Solving it subject to the initial condition $P_m(1,m)=1$ yields
\begin{equation}
\label{PdN:sol-m}
P_m(d,N)=(m-1)\frac{\Gamma(N-d)\,\Gamma(N-m+1)}{\Gamma(N-d-m+2)\,\Gamma(N)}
\end{equation}
\end{subequations}
The random quantity $d=d[v_m(N)]$ is non-self-averaging for every $m$. 

Alternatively, the results \eqref{PdN:sols} could be appreciated after realizing that the evolution of a degree of any vertex is equivalent to the P\'olya urn process  \cite{urn,urn-Mahmoud}. The P\'olya urn model proposed by Markov \cite{Markov17} and by Eggenberger and P\'olya \cite{Polya23} is the simplest, and best-understood urn model \cite{urn,urn-Mahmoud}. Urn models have been used by Huygens, Bernoulli, Laplace \cite{Laplace} and other founders of probability theory; according to \cite{urn}, traces of urn schemes appear already in the {\em Old Testament}. 

Known behaviors of the P\'olya urn model allow one to extract some leadership characteristics. Suppose we seek the probability $\mathcal{D}_2$ that the degree of the second vertex exceeds half the size of the hypergraph in the quickest possible way and then holds throughout the evolution:
\begin{equation}
\label{D2:def}
\mathcal{D}_2 = \frac{1}{2}\text{Prob}\{d[v_2(N)] > \lfloor N/2\rfloor \,| N\geq 3\}
\end{equation}
The factor $\frac{1}{2}$ accounts for the quickest path, viz. creating the hypergraph \eqref{123} with $d[v_2(3)] = 2>\lfloor 3/2\rfloor$. The probability in \eqref{D2:def} can be extracted from \cite{Tibor10} giving 
\begin{equation}
\label{D2}
\mathcal{D}_2 =\frac{1}{4}
\end{equation}

The virtue of $\mathcal{D}_2$ is in tractability and shedding light on a more natural quantity $\mathcal{S}_2$, the probability that the second vertex is the {\em strict} leader in degree (i.e., it has the second highest degree) for all sufficiently large sizes. We haven't computed $\mathcal{S}_2$, but the obvious inequality
\begin{equation}
\label{S2}
\mathcal{S}_2 > \mathcal{D}_2 =\frac{1}{4}
\end{equation}
tells us that with a positive probability the second vertex eventually becomes the persistent strict leader in degree.

Similarly, denote by $\mathcal{D}_m$ the probability that the degree of the $m^\text{th}$ vertex exceeds half the size of the hypergraph in the quickest possible way and then holds throughout the entire evolution. More precisely,
\begin{equation}
\label{Dm:def}
\mathcal{D}_m = \frac{\text{Prob}\{d[v_m(N)] > \lfloor N/2\rfloor \,| N\geq 2m-1\}}{(2m-2)!}
\end{equation}
with $1/(2m-2)!$ factor accounting for the quickest path, namely, creating the hypergraph
\begin{equation}
\label{quick-m}
\begin{split}
 & \{v_1\} \\
 & \{v_1,  v_2\}, ~ \cdots, ~\{v_1, v_m\}   \\
 & \{v_1, v_m, v_{m+1}\}, ~ \cdots, ~\{v_1, v_m, v_{m+1},\ldots,v_{2m-1}\}   
 \end{split}
\end{equation}
of the smallest size $N=2m-1$ when the inequality in \eqref{D2:def} is feasible, $d[v_m(2m-1)] = m>\lfloor (2m-1)/2\rfloor$. The probability $\mathcal{D}_m$ reads
\begin{equation}
\label{Dm}
\mathcal{D}_m =\frac{1}{(2m-2)!}\,\frac{\Gamma\big(m-\frac{1}{2}\big)}{\Gamma\big(\frac{1}{2}\big)\,\Gamma(m)} = \big[2^{m-1}\Gamma(m)\big]^{-2}
\end{equation}
The first formula for $\mathcal{D}_m$ is extracted from \cite{Tibor10}, and we reduced it to the second  by using the duplication formula \cite{Flajolet} for the gamma function.

A more natural quantity is again $\mathcal{S}_m$, the probability that the $m^\text{th}$ vertex is the strict leader in degree for all sufficiently large sizes. The obvious inequality 
\begin{equation}
\label{Sm}
\mathcal{S}_m > \mathcal{D}_m = \big[2^{m-1}\Gamma(m)\big]^{-2}
\end{equation}
tells us that with a positive probability, the $m^\text{th}$ vertex eventually becomes the persistent strict leader in degree. 

The monotonicity of the probabilities $\mathcal{S}_m$ is obvious: $\mathcal{S}_m>\mathcal{S}_{m+1}$ for all $m\geq 2$. Leading degrees can undergo a bit of leapfrogging. Conjecturally, the degree of one vertex eventually becomes the winner, namely the strict leader in degree for all sufficiently large sizes. If true, 
\begin{equation}
\label{norm-S}
\sum_{m\geq 2} \mathcal{S}_m = 1
\end{equation}
Replacing $\mathcal{S}_m$ in \eqref{norm} by the lower bound in \eqref{Sm} gives the lower bound for the sum, $I_0(1)-1=0.266\ldots$, where $I_0$ is the Bessel function. This lower bound for the sum \eqref{norm-S} is significantly lower than the exact value indicating that the lower bound \eqref{Sm} is rather weak.

\section{Leaves}
\label{sec:leaves}

The RRH model is too simple to be overlooked, and it appeared (albeit not called the RRH) in \cite{KR05,Vazquez22} and perhaps in other studies. More precisely, Ref.~\cite{Vazquez22} examines the one-parameter class of models: The RRH rule is applied with probability $1-\mu$, and with probability $\mu$, an edge is chosen randomly and duplicated. Thus after $N-1$ steps, there are $N$ edges, while the number of vertices is a random quantity concentrating around the average $\langle |V|\rangle = 1 +(N-1)(1-\mu)$. This class of models is a deformation of the RRH recovered at $\mu=0$. The influence of the deformation is minimal, e.g., the fractions $n_k=N_k(N)/N$ still follow the decay law \eqref{nk-sol} independently on $\mu$, albeit the stationarity is lost: When $0<\mu<1$, the decay law \eqref{nk-sol} is valid only when $N\to\infty$ and $k=O(1)$. Thus, the case of the RRH is particularly striking, and it well represents the behavior of the entire class of models \cite{Vazquez22}. 

\begin{figure}[t]
\includegraphics[width=0.44\textwidth]{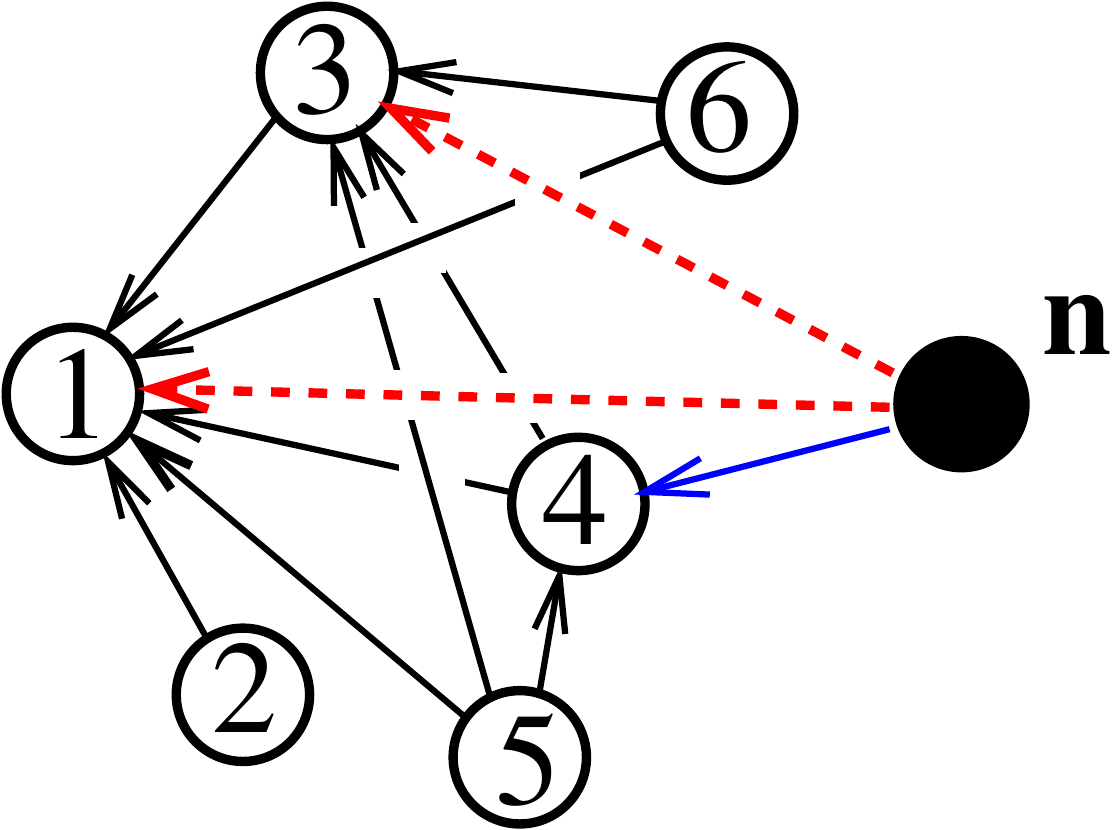}
\caption{Illustration of the copying mechanism. The time order of the nodes is indicated. The new node ${\bf n}$ chooses the target node randomly (node 4 in the present case) and attaches to it and its ancestors, nodes 1 and 3 in the present case.} 
\label{fig:C-color}
\end{figure}

The RRHs also appeared in an earlier work \cite{KR05} studying directed random graphs growing via a copying mechanism (CM). These graphs grow by adding nodes one by one. A newly introduced node randomly selects a target node and forms a direct link to it, as well as to all ancestors nodes of the target node (see Fig.~\ref{fig:C-color}). Identifying (i) nodes in a directed graph growing via CM with vertices in a hypergraph and (ii) each node and its ancestors with an edge, we establish the isomorphism between directed random graphs growing via CM and RRH. For the directed graph shown in Fig.~\ref{fig:C-color}, the corresponding hypergraph is
\begin{equation}
\label{CM-H}
\begin{split}
 & \{v_1\} \\
 & \{v_1,  v_2\}, ~\{v_1, v_3\}   \\
 & \{v_1, v_3,v_4\}, ~\{v_1, v_3,v_6\}  \\
 & \{v_1, v_3,v_4, v_5\}, ~\{v_1, v_3, v_4,v_7\}
 \end{split}
\end{equation}

Some results of Ref.~\cite{KR05} are identical to the results presented above, albeit the interpretation differs: The degree distribution for RRH was the in-degree distribution in Ref.~\cite{KR05}; the rank distribution for the RRH was the out-degree distribution in \cite{KR05}. The concepts of the degree distribution and rank distribution in the realm of the RRH seem more fundamental than the concepts of in- and out-degree distributions for graphs. Other characteristics could be more natural for graphs than for hypergraphs. For instance, leaves (nodes of degree one) are easily visible in connected graphs. For the graph in Fig.~\ref{fig:C-color}, node 2 is the only leaf. In RRH, a leaf is a vertex that belongs to a single edge, necessarily an edge of size two, since every edge contains the primordial vertex. The hypergraph shown in \eqref{E14} has two leaves, $v_8$ and $v_9$. 

Let us look at the number of leaves $\mathcal{L}_N$ in an RRH of size $N$. This number is a strongly fluctuating random quantity with average $\langle \mathcal{L}_N\rangle = 1$ for $N\geq 3$ when the number of leaves is a random quantity. Indeed
\begin{equation}
\mathcal{L}_3=
\begin{cases}
0 & \text{prob}~~\frac{1}{2}\\
2 & \text{prob}~~\frac{1}{2}
\end{cases}
\end{equation}
from which $\langle \mathcal{L}_3\rangle = 1$. Further 
\begin{equation}
\label{LNN}
\mathcal{L}_{N+1}=
\begin{cases}
\mathcal{L}_N-1 & \text{prob}~~\frac{\mathcal{L}_N}{N}\\
\mathcal{L}_N+1 & \text{prob}~~\frac{1}{N}\\
\mathcal{L}_N      & \text{otherwise}
\end{cases}
\end{equation}
for $N\geq 3$. Averaging \eqref{LNN} yields
\begin{equation}
\langle \mathcal{L}_{N+1}\rangle= \langle \mathcal{L}_N\rangle+ \frac{1-\langle \mathcal{L}_N\rangle}{N}
\end{equation}
Starting with $\langle \mathcal{L}_3\rangle = 1$ we iterate and obtain $\langle \mathcal{L}_N\rangle = 1$ for all $N\geq 3$. 

The full probability distribution 
\begin{equation}
\ell_k(N) = \text{Prob}[\mathcal{L}_N=k]
\end{equation}
satisfies an equation
\begin{eqnarray}
\label{LkN}
\ell_k(N+1) &=& \left[1-\frac{k+1}{N}\right]\ell_k(N) \nonumber \\
&+&\frac{1}{N}\,\ell_{k-1}(N) + \frac{k+1}{N}\,\ell_{k+1}(N)
\end{eqnarray}
which is derived similarly to Eq.~\eqref{LNN}. In the $N\to\infty$ limit, the probabilities $\ell_k(N)$ saturate for any fixed $k$:
\begin{equation}
\ell_k \equiv \lim_{N\to\infty} \ell_k(N)
\end{equation}
Using \eqref{LkN} we deduce the recurrence 
\begin{equation}
\label{Lk-rec}
\ell_{k+1}=\ell_k-(k+1)^{-1}\ell_{k-1}
\end{equation}
for $k\geq 1$. Solving \eqref{Lk-rec} starting with $\ell_0=\ell_1$ which also follows from \eqref{LkN} we obtain
\begin{equation}
\label{Lk-sol}
\ell_k = \frac{e^{-1}}{k!}
\end{equation}
with amplitude fixed by normalization: $\sum_{k\geq 0}\ell_k = 1$. 

\section{Redirection Mechanism}
\label{sec:RD}

Deformations of the RRT model via different types of redirection (see \cite{Kleinberg,KR01,KR02-fluct,Slava05a,KR05,BK10,KR-copy,KR14,Bertoin15,KR17,Biham19a, Biham19b,Sumpter22} and references therein), have significantly improved our understanding of growing random networks. Natural deformations of the RRH model may also lead to models of growing random hypergraphs exhibiting intriguing behaviors.

Let us deform the RRH using the simplest implementation of the redirection mechanism.  At each step, we add a vertex $v$ and an edge formed from a randomly chosen existing edge $e$. This new edge is $e\cup {v}$ with probability $1-r$, while with probability $r$, we redirect from $e$ to its maternal edge $e'$, so the new edge is $e'\cup {v}$. If $e=\{v_{i_1}, \ldots,v_{i_n}\}$ with $1=i_1<\ldots<i_n$, we have $e'=\{v_{i_1}, \ldots,v_{i_{n-1}}\}$. Thus the added edge is
\begin{equation}
\label{RD}
\begin{cases}
\{v_{i_1}, \ldots,v_{i_{n-1}},v_{i_n}, v\} & \text{prob}\quad 1-r \\
\{v_{i_1}, \ldots,v_{i_{n-1}},v\}              & \text{prob}\quad r
\end{cases}
\end{equation}
The primordial vertex has no maternal edge, so if the primordial vertex is chosen,  $e=\{v_1\}$,   the added edge is always $e\cup {v}=\{v_1,v\}$. 

Consider vertices of degree one. The governing relation \eqref{N1N} generalizes to
\begin{equation}
\label{N1N-r}
\mathcal{N}_1(N+1) = 
\begin{cases}
\mathcal{N}_1+1 & \text{prob} ~~1-(1-r)\frac{\mathcal{N}_1}{N} \\
\mathcal{N}_1    & \text{prob} ~~ (1-r)\frac{\mathcal{N}_1}{N}
\end{cases}
\end{equation}
for $N\geq 2$. Averaging we obtain
\begin{equation}
\label{N1-r}
N_1(N+1)=\left(1-\frac{1-r}{N}\right)N_1+1
\end{equation}
for $N\geq 2$. Solving \eqref{N1-r} subject to the boundary condition $N_1(2)=1$ we obtain
\begin{equation}
\label{N1-sol-r}
N_1 = \frac{N}{2-r}-\frac{1}{(2-r)\Gamma(r)}\,\frac{\Gamma(N-1+r)}{\Gamma(N)}
\end{equation}
for $N\geq 2$. In contrast to the RRH model with $r=0$ when $N_1$ was strictly linear in $N$, there is a sub-leading term in \eqref{N1-sol-r} when $0<r<1$ that vanishes as $N^{-(1-r)}$. 

For $k\geq 2$, the random quantity $\mathcal{N}_k$ evolves according to stochastic rule 
\begin{equation*}
\label{NkN-r}
\mathcal{N}_k(N+1) =
\begin{cases}
\mathcal{N}_k+1  &\text{prob} ~~ \frac{(k-1-r)\mathcal{N}_{k-1}}{N}\\
\mathcal{N}_k-1  &\text{prob} ~~ \frac{(k-r)\mathcal{N}_{k}}{N} \\
\mathcal{N}_k     &\text{prob} ~~  1-\frac{(k-1-r)\mathcal{N}_{k-1}+(k-r) \mathcal{N}_k}{N}
\end{cases}
\end{equation*}
Averaging this equation and taking the $N\to\infty$ limit we find that fractions $n_k=N_k(N)/N$ satisfy
\begin{equation}
\label{nkr}
(k-r+1)n_k = (k-r-1)n_{k-1}
\end{equation}
Using $n_1=(2-r)^{-1}$ following from \eqref{N1-sol-r} as the initial condition, we solve \eqref{nkr} recurrently and find
\begin{equation}
\label{nkr:sol}
n_k = \frac{1-r}{(k-r+1)(k-r)}
\end{equation}

We see that for the one-parameter class of models \eqref{RD}, the redirection parameter $r$ only quantitatively affects the behavior of the degree distribution. The influence on the degree distribution is more significant than in the one-parameter class of models introduced in  \cite{Vazquez22}, but still just quantitative. However, the degree distribution in growing hypergraphs appears to be a very robust characteristic hardly sensitive to the evolution rules. Despite the name, the degree distribution in the hypergraphs, Eq.~\eqref{Nk:def} resembles the in-component size distribution in trees rather than the degree distribution. Intriguingly, in the one-parameter class of models of trees growing via the redirection mechanism \cite{KR01}, the in-component size distribution is given {\em exactly} by \eqref{nkr:sol}. The degree distribution in such trees is strongly affected by the parameter $r$: It decays algebraically as $k^{-1-1/r}$ when $0<r<1$, while for the RRTs ($r=0$) the decay is exponential, $n_k=2^{-k}$. 

The influence of the redirection parameter $r$ on the rank distribution is more substantial. Consider the vertices of rank two. For the RRHs, Eq.~\eqref{R2-sol} gives the average number of vertices of rank two. The growth with the size of the hypergraphs is logarithmic, Eq.~\eqref{R2-asymp}. When redirection can occur, $0<r<1$, the quantity $R_2(N)$ grows with $N$ algebraically, as we now demonstrate. 

To establish the growth law we write the stochastic equation for the number of vertices of rank two:
\begin{equation}
\label{R2N-r}
\mathcal{R}_2(N+1) = 
\begin{cases}
\mathcal{R}_2+1 & \text{prob} ~~\frac{1+r \mathcal{R}_2}{N}\\
\mathcal{R}_2    & \text{prob} ~~ 1-\frac{1+r \mathcal{R}_2}{N}
\end{cases}
\end{equation}
The average $R_2=\langle \mathcal{R}_2\rangle$ satisfies the recurrence
\begin{equation}
\label{R2-rec-r}
R_2(N+1) = \left(1+\frac{r}{N}\right)R_2 + \frac{1}{N}
\end{equation}
from which
\begin{equation}
\label{R2-sol-r}
R_2(N) = \frac{1}{r}\left[\frac{\Gamma(N+r)}{\Gamma(1+r)\,\Gamma(N)} -1\right]
\end{equation}
Using the large $N$ asymptotic, $\Gamma(N+r)/\Gamma(N)\to N^r$, for the ratio of gamma functions \cite{Knuth}, we deduce an 
algebraic growth law for the average number of vertices of rank two: $R_2\simeq \frac{N^r}{r\,\Gamma(1+r)}$ when $0<r<1$. 

Similarly $\mathcal{R}_k$ with $k\geq 3$ satisfy 
\begin{equation}
\label{RkN-r}
\mathcal{R}_k(N+1) = 
\begin{cases}
\mathcal{R}_k+1 & \text{prob} ~~\frac{(1-r)\mathcal{R}_{k-1}+r\mathcal{R}_k}{N}\\
\mathcal{R}_k    & \text{prob} ~~ 1- \frac{(1-r)\mathcal{R}_{k-1}+r\mathcal{R}_k}{N}
\end{cases}
\end{equation}
from which we deduce the recurrence for the averages:
\begin{equation}
\label{Rk-rec-r}
R_k(N+1) = \left(1+\frac{r}{N}\right)R_k + \frac{1-r}{N}\,R_{k-1}
\end{equation}

Using the exact solution \eqref{R2-sol-r} one finds $R_3$, then $R_4$, etc. These explicit exact results are cumbersome, so we merely give the leading asymptotic:
\begin{equation}
\label{Rk-asymp-r}
R_{k+2} \simeq \frac{N^r}{r\,\Gamma(1+r)}\,\frac{[(1-r)\ln N]^k}{k!}
\end{equation}
This asymptotic can be derived by replacing the recurrence \eqref{Rk-rec-r} by a system of differential equations and using the aforementioned leading behavior of $R_2$.

\section{Discussion}
\label{sec:concl}

We have investigated random recursive hypergraphs (RRHs) built via simple growth rules. The RRH model is parameter-free, so it can be considered a null model of growing random hypergraphs. Several characteristics of RRHs admit neat expressions via beautiful special numbers (harmonic numbers, Bernoulli numbers, Eulerian numbers, and Stirling numbers of the first kind) and are valid for arbitrary $N$. These exact results depend on the initial condition. Many basic random quantities characterizing the RRHs are asymptotically self-averaging, so their leading asymptotic behaviors are independent on the initial condition. Some random quantities are non-self-averaging, e.g., the leadership characteristics of the RRHs, so their  asymptotic behaviors fluctuate from realization to realization.

We have also briefly looked at a one-parameter class of models, a simple deformation of the RRH model. Namely, with a certain probability $r$, a parameter of the model, the redirection from a selected edge to its maternal edge is allowed. A detailed analysis of this class of models is a natural direction for future research. In the realm of random graphs, the redirection mechanism \cite{Kleinberg,KR01} generates preferential attachment. When redirection not only to the closest ancestor is allowed, the formation of hubs becomes feasible (see \cite{BK10}). It would be interesting to investigate hypergraphs built via such generalized redirection. For undirected trees, a particularly striking behavior was found in the case of isotropic complete redirection \cite{KR17}. A hypergraph version corresponds to $r=1$ and uniform choice among all neighbors of the initially chosen edge $e$. 

Complex networks growing via choice-driven rules exhibit phase transitions and other unexpected behaviors \cite{Raissa07,Mahmoud10,KR14-choice,Malyshkin14,Jordan16}. The same could happen for choice-driven growing hypergraphs. A simple implementation of choice relies on provisionally selecting two edges at random, say $e_1$ and $e_2$; choosing one of them, say $e_i$, according to some rule; and adding a new vertex $v$ together with edge $e_i\cup {v}$. For instance, the smaller edge is always chosen, so if $|e_1|<|e_2|$ the new edge is $e_1\cup {v}$; if $|e_1|=|e_2|$, the new edge is $e_i\cup {v}$ where $i=1$ or $i=2$ with equal probabilities. 

An intriguing direction of future research concerns growing densifying hypergraphs for which the number of edges grows qualitatively faster than the number of vertices.  A hypergraph version of the growing network model \cite{KR-copy} is expected to lead to growing densifying hypergraphs. As for RRHs, one randomly chooses an edge $e$ and adds an edge $e\cup\{v\}$; in addition, for each ancestor $e'\subset e$, an edge $e'\cup\{v\}$ is added with probability $p$. Numerous biological and technological networks are densifying \cite{brain05,brain09,brain11,Kleinberg07}. Densifying hypergraphs are also widespread, so simple models of growing densifying hypergraphs may shed light on the properties of such objects.  

Growing graphs are sparse if the ratio of the number of edges to the number of nodes remains finite; if $|E|/N$ diverges as $N\to \infty$, growing graphs are densifying. For connected graphs, $|E|/N$ varies from one for trees (more precisely, $|E|=N-1$ for trees) to $\frac{N-1}{2}$ characterizing the complete graph of size $N$. For hypergraphs, the maximal number of edges is $2^{N}-1$. Thus, a potent sparse vs. dense dichotomy is not necessarily the same \footnote{For hypergraphs, the ratio of logarithms  could be a better quantifier. Log-densifying hypergraphs are those for which the ratio $\ln |E|/\ln N$ diverges in the $N\to \infty$ limit.} for graphs and hypergraphs.

\bigskip
\noindent{\bf Acknowledgments.}
I am grateful to G. Bianconi, H. Hartle, D. Krioukov, S. Redner, and J. Stepanyants for the discussions.

\bibliography{references-nets}

\end{document}